\documentclass[a4paper]{amsart}

\usepackage{amsmath, tabularx}
\usepackage{amsthm}
\usepackage{amssymb}
\usepackage{amsfonts}
\usepackage{paralist}
\usepackage{aliascnt}
\usepackage{amsrefs}
\usepackage{amscd}
\usepackage{blkarray}

\usepackage{setspace}
\usepackage[colorlinks=true,linkcolor=blue,urlcolor=blue]{hyperref}
\usepackage{tikz}
\usetikzlibrary{matrix}
\usetikzlibrary{arrows}

\newcommand{\kk}{\normalfont\mathbb{K}}

\newcommand{\reg}{\normalfont\text{reg}}

\newcommand{\Tor}{\normalfont\text{Tor}}
\newcommand{\Ext}{\normalfont\text{Ext}}

\newcommand{\Ker}{\normalfont\text{Ker}}

\newcommand{\Rees}{\mathcal{R}}

\newcommand{\EEQ}{\mathcal{K}}

\newcommand{\iniTerm}{\normalfont\text{in}}
\newcommand{\bideg}{\normalfont\text{bideg}}

\newcommand{\match}{\mathcal{M}}
\newcommand{\GB}{\mathcal{G}_{<^\match}(\EEQ)}

\newtheorem{theorem}{Theorem}[section]

\newaliascnt{corollary}{theorem}
\newtheorem{corollary}[corollary]{Corollary}
\aliascntresetthe{corollary}

\newaliascnt{lemma}{theorem}
\newtheorem{lemma}[lemma]{Lemma}
\aliascntresetthe{lemma}

\newaliascnt{conjecture}{theorem}
\newtheorem{conjecture}[conjecture]{Conjecture}
\aliascntresetthe{conjecture}

\newaliascnt{proposition}{theorem}
\newtheorem{proposition}[proposition]{Proposition}
\aliascntresetthe{proposition}

\newaliascnt{definition}{theorem}
\newtheorem{definition}[definition]{Definition}
\aliascntresetthe{definition}

\newaliascnt{notation}{theorem}
\newtheorem{notation}[notation]{Notation}
\aliascntresetthe{notation}

\newaliascnt{example}{theorem}
\newtheorem{example}[example]{Example}
\aliascntresetthe{example}

\newaliascnt{remark}{theorem}
\newtheorem{remark}[remark]{Remark}
\aliascntresetthe{remark}

\newaliascnt{problem}{theorem}

\aliascntresetthe{problem}

\newaliascnt{construction}{theorem}

\aliascntresetthe{construction}

\newaliascnt{defprop}{theorem}

\aliascntresetthe{defprop}

\def\equationautorefname~#1\null{(#1)\null}
\def\sectionautorefname~#1\null{Appendix\null}

\begin{document}

\title[Regularity and Gr\"obner bases of $\Rees(I)$ for bipartite graphs]{Regularity and Gr\"obner bases of the Rees algebra of edge ideals of bipartite graphs}
\author{Yairon Cid-Ruiz}
\address{Department de Matem\`{a}tiques i Inform\`{a}tica, Facultat de Matem\`{a}tiques i Inform\`{a}tica, Universitat de Barcelona, Gran Via de les Corts Catalanes, 585; 08007 Barcelona, Spain.}
\email{ycid@ub.edu}
\urladdr{http://www.ub.edu/arcades/ycid.html}




\thanks{The author was funded by the European Union's Horizon 2020 research and innovation programme under the Marie Sk\l{}odowska-Curie grant agreement No. 675789.
The author acknowledges financial support from the Spanish Ministry of
Economy and Competitiveness, through the ``Mar\'ia de Maeztu'' Programme for Units of Excellence in R\&D (MDM-2014-0445).}

\begin{abstract}
	Let $G$ be a bipartite graph and $I=I(G)$ be its edge ideal. 
	The aim of this note is to investigate different aspects of the Rees algebra $\Rees(I)$ of $I$. 
	We compute its regularity and the universal Gr\"obner basis of its defining equations; interestingly, both of them are described in terms of the combinatorics of $G$.

	We apply these ideas  to study the regularity of the powers of $I$.
	For any $s \ge \text{match}(G)+\lvert E(G) \rvert +1$ we prove that $\reg(I^{s+1})=\reg(I^s)+2$.
\end{abstract}
\maketitle
\vspace*{-.9cm}

\section{Introduction}

Let $G=\left(V(G), E(G)\right)$ be a bipartite graph on the vertex set $V(G)=X \cup Y$ with bipartition $X=\{x_1, \ldots, x_n \}$ and $Y=\{y_1,\ldots,y_m\}$. 
Let $\kk$ be a field  and let $R$ be the polynomial ring $R = \kk[x_1,\ldots ,x_n,y_1,\ldots,y_m]$.
The edge ideal $I=I(G)$, associated to $G$, is the ideal of $R$ generated by the set of monomials $x_iy_j$ such that $x_i$ is adjacent to $y_j$.

One can find a vast literature on the Rees algebra of edge ideals of bipartite graphs (see \cite{MONOMIAL_ALGEBRAS}, \cite{SIMIS_VASC_VILLARREAL_IDEAL_GRAPH}, \cite{VILLARREAL_BIPARTITE}, \cite{VILLARREAL_COMPLETE_BIPART}, \cite{VILLARREAL_EDGE}, \cite{REES_ALGEBRAS_POLYH_CONES}, \cite{LOUIZA_KUEI}), nevertheless, in this note we study several properties that might have been overlooked.
From a computational point of view we first focus on the universal Gr\"obner basis of its defining equations, and from a more algebraic standpoint we focus on its total and partial regularities as a bigraded algebra. 
Applying these ideas, we give an estimation of when $\reg(I^s)$ starts to be a linear function and we find upper bounds for the regularity of the powers of $I$.

Let $\Rees(I)=\bigoplus_{i=0}^{\infty} I^it^i \subset R[t]$ be the Rees algebra of the edge ideal $I$.
Let $f_1,\ldots,f_q$ be the square free monomials of degree two generating $I$.
We can see $\Rees(I)$ as a quotient of the polynomial ring $S=R[T_1,\ldots,T_q]$ via the map
\begin{align}
\label{def_Rees}
\begin{split}
S&=\kk[x_1,\ldots,x_n,y_1\ldots,y_m,T_1,\ldots,T_q] \xrightarrow{\psi} \Rees(I) \subset R[t],\\ \quad
&\psi(x_i) = x_i,\quad \psi(y_i) = y_i, \quad \psi(T_i) = f_it.
\end{split}
\end{align}
Then the presentation of $\Rees(I)$ is given by $S/\EEQ$ where $\EEQ = \Ker(\psi)$.
We give a bigraded structure to $S=\kk[x_1,\ldots,x_n,y_1,\ldots,y_m] \otimes_{\kk} \kk[T_1,\ldots,T_q]$, where $\bideg(x_i)=\bideg(y_i)=(1,0)$ and $\bideg(T_i)=(0, 1)$.  
The map $\psi$ from \autoref{def_Rees} becomes bihomogeneous
when we declare $\bideg(t)=(-2,1)$, then we have that $S/\EEQ$ and $\EEQ$ have natural bigraded structures as $S$-modules.

The universal Gr\"obner basis of the ideal $\EEQ$ is defined as the union of all the reduced Gr\"obner bases $\mathcal{G}_<$ of the ideal $\EEQ$ as $<$ runs over all possible monomial orders (see \cite{STURMFELS_MONOG}).
In our first main result we compute the universal Gr\"obner basis of the defining equations $\EEQ$ of the Rees algebra $\Rees(I)$.

\begin{theorem}[\autoref{thm_univ_grobner_basis}]
	\label{thmA}
	Let $G$ be a bipartite graph and $\EEQ$ be the defining equations of the Rees algebra $\Rees(I(G))$. 
	The universal Gr\"obner basis $\mathcal{U}$ of $\EEQ$ is given by 
	\begin{align*}
	\mathcal{U} &=  \{ T_w \mid w \text{ is an even cycle}\} \\
	&\cup  \{ v_0T_{w^+}-v_aT_{w^-} \mid w=(v_0,\ldots,v_a) \text{ is an even path} \}\\
	&\cup \{u_0u_aT_{(w_1,w_2)^+}- v_0v_bT_{(w_1,w_2)^-} \mid w_1=(u_0,\ldots,u_a) \text{ and }\\ &\qquad\qquad\qquad\qquad\qquad\qquad w_2=(v_0,\ldots,v_b) \text{ are disjoint odd paths}\}.	
	\end{align*}
\end{theorem}

From \cite[Theorem 3.1, Proposition 3.1]{VILLARREAL_EDGE} we have a precise description of $\EEQ$ given by the syzygies of $I$ and  the set even of closed walks in the graph $G$.
The algebra $\Rees(I)$, as a bigraded $S$-module, has a minimal bigraded free resolution 
\begin{equation}
\label{minimal_free_res_Rees}
0 \longrightarrow F_p \longrightarrow \cdots \longrightarrow F_1 \longrightarrow F_0 \longrightarrow \Rees(I) \longrightarrow 0,
\end{equation}
where $F_i=\oplus_j S(-a_{ij},-b_{ij})$. 
In the same way as in \cite{ROMER_BIGRAD},  we can define the $xy$-regularity of $\Rees(I)$ by the integer
$$
\reg_{xy}(\Rees(I)) = \max_{i,j} \{a_{ij} - i\},
$$
or equivalently by
$$
\reg_{xy}(\Rees(I)) = \max\{ a \in \mathbb{Z} \mid \beta_{i, (a+i, b)}^S (\Rees(I)) \neq 0 \text{ for some } i,b \in \mathbb{Z} \},
$$
where $\beta_{i, (a, b)}^S (\Rees(I)) = \dim_{\kk}({\Tor_i^S(\Rees(I), \kk)}_{(a,b)})$.

Similarly, we can define the $T$-regularity 
$$
\reg_{T}(\Rees(I)) = \max_{i,j} \{b_{ij} - i\}
$$
and the total regularity 
$$
\reg(\Rees(I)) = \max_{i,j} \{a_{ij}+b_{ij} - i\}.
$$

Our second main result is computing the total regularity and giving upper bounds for both partial regularities. 

\begin{theorem}[\autoref{reg_Rees}]
	\label{thmB}
	Let $G$ be a bipartite graph. Then we have:
	\begin{enumerate}[(i)]
		\item \; $\reg(\Rees(I(G))) = \text{match}(G)$,
		\item \; $\reg_{xy}(\Rees(I(G))) \le \text{match}(G) -1$,
		\item \;$\reg_{T}(\Rees(I(G))) \le \text{match}(G)$,
	\end{enumerate}
 	where $\normalfont\text{match}(G)$ denotes the matching number of $G$.
\end{theorem}

Finally, we apply these results in order to study the regularity of the powers of the edge ideal $I=I(G)$. 

It is a famous result (for a general ideal in a polynomial ring) the asymptotic linearity of $\reg(I^s)$ for $s\gg0$ (see \cite{CUTKOSKY_HERZOG_TRUN_LINEARITY_REG} and \cite{VIJAY_REG_LINEAR}). 
However, the exact form of this linear function and the exact point where $\reg(I^s)$ starts to be linear, is a problem that continues wide open even in the case of monomial ideals.

In recent years, a number of researchers have focused on computing the regularity of powers of edge ideals and on relating these values to combinatorial invariants of the graph (see e.g. \cite{BBH}, \cite{REG_UNICYCLIC_GRAPHS}, \cite{Three_Bipartite_Graphs}, \cite{REG_POWERS_EDGE_BANERJEE}, \cite{FORESTS_CYCLES}, \cite{REG_POWERS_BIPART}).
Most of the upper bounds given in these papers use the concept of even-connection introduced in \cite{REG_POWERS_EDGE_BANERJEE}.
Actually, using this idea as a central tool, in \cite{REG_POWERS_BIPART} it was proved the upper bound 
$$
\reg(I^s) \le 2s + \text{co-chord}(G)-1
$$
for any bipartite graph $G$, where $\text{co-chord}(G)$ represents the co-chordal number of $G$ (see \cite[Definition 3.1]{REG_POWERS_BIPART}).

As a consequence of our study of the Rees algebra $\Rees(I)$,  we make an estimation of when $\reg(I^s)$ starts to be a linear function, and we obtain the weaker upper bounds for the regularity of the powers of $I$ (see \autoref{comparison}, \autoref{cor_upper_bound_all_powers}, \autoref{cor_upper_bound_bipart_order}). 
Perhaps, this could give new tools and fresh ideas to pursue the stronger upper bound
\begin{equation}
	\label{upper_bound_Ha}
	\reg(I^s) \le 2s + \reg(I) -2,
\end{equation}
that has been conjectured by  Alilooee, Banerjee, Beyarslan and H\`a (\cite[Conjecture 7.11]{BBH}).


Using the upper bound for the partial $T$-regularity of $\Rees(I)$, we can get the following estimation.

\begin{corollary}[\autoref{cor_eventual_linearity}]
	\label{thmD}
	Let $G$ be a bipartite graph. Then, for all $s \ge {\normalfont\text{match}}(G)+\lvert E(G) \rvert+1$ we have
		$$
		\reg({I(G)}^{s+1})=\reg({I(G)}^s)+2.
	    $$	
\end{corollary}


The basic outline of this note is as follows. 

In \hyperref[section2]{Section 2}, we compute the universal Gr\"obner basis of $\EEQ$ (\autoref{thmA}).
In \hyperref[section3]{Section 3}, we consider a specific monomial order that allows us to get upper bounds for the $xy$-regularity of $\Rees(I)$.
In \hyperref[section4]{Section 4} we exploit the canonical module of $\Rees(I)$ in order to prove \autoref{thmB} and  \autoref{thmD}. 
Finally, in \hyperref[section5]{Section 5} we give some general ideas about the conjectured upper bound  \autoref{upper_bound_Ha}.

\section{The universal Gr\"obner basis of $\EEQ$}\label{section2}

In this section we will give an explicit description of the universal Gr\"obner basis $\mathcal{U}$ of $\EEQ$.
Our approach is the following, 
first we compute the set of circuits of the incidence matrix of the cone graph, and then we translate this set of circuits into a description of $\mathcal{U}$.

The following notation will be assumed in most of this note.
\begin{notation}
	\label{general_notation}
		Let $G$ be a bipartite graph with bipartition $X=\{x_1,\ldots,x_n\}$ and $Y=\{y_1,\ldots,y_m\}$, and $R$ be the polynomial ring $R=\kk[x_1,\ldots,x_n,y_1,\ldots,y_m]$.
	Let $I$ be the edge ideal $I(G)=\left(f_1,\ldots,f_q\right)$ of $G$.
	We consider the Rees algebra $\Rees(I)$ as a quotient of $S=R[T_1,\ldots,T_q]$ by using  \autoref{def_Rees}. 
	Let $\EEQ$ be the defining equations of the Rees algebra $\Rees(I)$.
\end{notation}


Let $A=(a_{i,j}) \in \mathbb{R}^{n+m,q}$ be the incidence matrix of the graph $G$. 
Then we construct the matrix $M$ of the following form 
\begin{equation}
\label{matrix_def_Rees}
M = \left( 
\begin{array}{cccccc}
a_{1,1} & \ldots & a_{1,q} & e_1 & \ldots & e_{n+m} \\
\vdots & \vdots & \vdots \\
a_{n+m,1} & \ldots & a_{n+m,q} \\
1 & \ldots & 1
\end{array}
\right),
\end{equation}
where $e_1, \ldots e_{n+m}$ are the first $n+m$ unit vectors in $\mathbb{R}^{n+m+1}$ (see \cite[Section 3]{VILLARREAL_BIPARTITE} for more details).
This matrix corresponds to the presentation of $\Rees(I)$ given in \autoref{def_Rees}.
For any vector $\beta \in \mathbb{Z}^{n+m+q}$ with nonnegative coordinates we shall use the notation
$$
{\mathbf{xyT}}^\beta = x_1^{\beta_{q+1}}\cdots x_n^{\beta_{q+n}} y_1^{\beta_{q+n+1}}\cdots y_m^{\beta_{q+n+m}} T_1^{\beta_1} \cdots T_q^{\beta_q}.
$$

A given vector $\alpha \in \Ker(M) \cap \mathbb{Z}^{n+m+q}$, can be written as $\alpha = \alpha^+ - \alpha^-$ where $\alpha^+$ and $\alpha^-$ are nonnegative and have disjoint support.

\begin{definition}[\cite{STURMFELS_MONOG}]
	A vector $\alpha \in \Ker(M) \cap \mathbb{Z}^{n+m+q}$
	is called a circuit if it has minimal support  $\text{supp}(\alpha)$ with respect to inclusion and its coordinates are relatively prime.	
\end{definition}

\begin{notation}
	\label{notation_trans}
	Given a  walk $w=\{v_0,\ldots,v_a\}$, each edge $\{v_{j-1},v_j\}$  corresponds to a variable $T_{i_j}$, and we set $T_{w^+}=\prod_{j \text{ is even}} T_{i_j}$ and $T_{w^-}=\prod_{j \text{ is odd}} T_{i_j}$ (in case $a=1$ we make $T_{w^+}=1$). 
	We adopt the following notations:
	\begin{enumerate}[(i)]
		\item Let $w = \{v_0, \ldots, v_a=v_0 \}$ be an even cycle in $G$. 
		Then by $T_w$ we will denote the binomial $T_{w^+}-T_{w^-} \in \EEQ$.
		\item Let $w=\{v_0, \ldots, v_a\}$ be an even path in $G$, since $G$ is bipartite then both endpoints of $w$ belong to the same side of the bipartition, i.e. either $v_0=x_i, v_a=x_j$ or $v_0=y_i, v_a=y_j$. 
		Then the path $w$ determines the binomial $$v_0T_{w^+}-v_aT_{w^-} \in \EEQ.$$ 
		\item Let $w_1=\{u_0,\ldots,u_a\}$, $w_2=\{v_0,\ldots,v_b\}$ be two disjoint odd paths, then the endpoints of $w_1$ and $w_2$ belong to different sides of the bipartition.
		Let $T_{(w_1,w_2)^+}=T_{w_1^+}T_{w_2^-}$ and $T_{(w_1,w_2)^-}=T_{w_1^-}T_{w_2^+}$, then $w_1$ and $w_2$ determine the binomial $$
		u_0u_aT_{(w_1,w_2)^+}- v_0v_bT_{(w_1,w_2)^-} \in \EEQ.$$
	\end{enumerate}
\end{notation}

\begin{example}
	\label{examp_graph_two_odd_cycles}
	In the bipartite graph shown below
	\begin{center}
		\definecolor{qqqqff}{rgb}{0.,0.,1.}
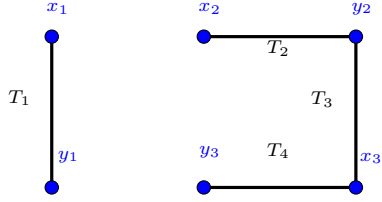
	\begin{tikzpicture}[line cap=round,line join=round,>=triangle 45,x=1.0cm,y=1.0cm]
	\clip(1.,0.5) rectangle (6.7,3.6);
	\draw [line width=1.2pt] (2.,3.)-- (2.,1.);
	\draw [line width=1.2pt] (4.,3.)-- (6.,3.);
	\draw [line width=1.2pt] (6.,3.)-- (6.,1.);
	\draw [line width=1.2pt] (6.,1.)-- (4.,1.);
	\begin{scriptsize}
	\draw [fill=qqqqff] (2.,3.) circle (2.5pt);
	\draw[color=qqqqff] (2.08,3.37) node {$x_1$};
	\draw [fill=qqqqff] (2.,1.) circle (2.5pt);
	\draw[color=qqqqff] (2.22,1.41) node {$y_1$};
	\draw[color=black] (1.58,2.19) node {$T_1$};
	\draw [fill=qqqqff] (4.,3.) circle (2.5pt);
	\draw[color=qqqqff] (4.08,3.37) node {$x_2$};
	\draw [fill=qqqqff] (6.,3.) circle (2.5pt);
	\draw[color=qqqqff] (6.08,3.37) node {$y_2$};
	\draw[color=black] (4.98,2.85) node {$T_2$};
	\draw [fill=qqqqff] (6.,1.) circle (2.5pt);
	\draw[color=qqqqff] (6.2,1.39) node {$x_3$};
	\draw[color=black] (5.56,2.17) node {$T_3$};
	\draw [fill=qqqqff] (4.,1.) circle (2.5pt);
	\draw[color=qqqqff] (4.08,1.47) node {$y_3$};
	\draw[color=black] (4.98,1.49) node {$T_4$};
	\end{scriptsize}
	\end{tikzpicture}
	\end{center}
	we have that the odd paths $w_1=(x_1,y_1)$ and $w_2=(x_2,y_2,x_3,y_3)$ determine the binomial
	$
	x_1y_1T_2T_4-x_2y_3T_1T_3.
	$
\end{example}

Let $\mathcal{U}$ be the universal G\"obner basis of $\EEQ$. 
In general we have that the set of circuits is contained in $\mathcal{U}$ (\cite[Proposition 4.11]{STURMFELS_MONOG}). 
But from the fact that $M$ is totally unimodular (\cite[Theorem 3.1]{VILLARREAL_BIPARTITE}), we can use \cite[Proposition 8.11]{STURMFELS_MONOG} and obtain the equality
$$
\mathcal{U} = \{ \mathbf{xyT}^{\alpha^+}
- \mathbf{xyT}^{\alpha^-} \mid \alpha \text{ is a circuit of } M\}.
$$


Therefore we shall focus on determining the circuits of $M$, and for this we will need to introduce the concept of the cone graph $C(G)$. 
The vertex set of the graph $C(G)$ is obtained by adding a new vertex $z$ to $G$, and its edge set consists of the edges in $E(G)$ together with the edges $\{x_1,z\},\ldots,\{x_n,z\}, \{y_1,z\}, \ldots, \{y_m,z\}$.

\begin{theorem}
	\label{thm_univ_grobner_basis}
	Let $G$ be a bipartite graph and $I=I(G)$ be its edge ideal.
	The universal Gr\"obner basis $\mathcal{U}$ of $\EEQ$ is given by 
	\begin{align*}
		\mathcal{U} &=  \{ T_w \mid w \text{ is an even cycle}\} \\
		&\cup  \{ v_0T_{w^+}-v_aT_{w^-} \mid w=(v_0,\ldots,v_a) \text{ is an even path} \}\\
		&\cup \{u_0u_aT_{(w_1,w_2)^+}- v_0v_bT_{(w_1,w_2)^-} \mid w_1=(u_0,\ldots,u_a) \text{ and }\\ &\qquad\qquad\qquad\qquad\qquad\qquad w_2=(v_0,\ldots,v_b) \text{ are disjoint odd paths}\}.	
	\end{align*}
	\begin{proof}
		Let $\kk[C(G)]$ be the monomial subring of the graph $C(G)$, which is generated by the monomials 
		$$
		\kk[C(G)] = \kk\big[ \{x_iy_j \mid \{x_i,y_j\} \in E(G) \} \cup \{x_iz \mid i=1,\ldots,n\} \cup \{y_iz \mid i=1,\ldots,m\}\big].
		$$
		
		As we did for the Rees algebra $\Rees(I)$, we can define a similar surjective homomorphism 
		\begin{align*}
			&\pi: S \longrightarrow \kk[C(G)] \subset R[z], \\
			\pi(x_i) &= x_iz,\quad \pi(y_i) = y_iz, \quad \pi(T_i) = f_i.\
		\end{align*}
		We have a natural isomorphism between $\Rees(I)$ and $\kk[C(G)]$ \cite[Excercise 7.3.3]{VASCONCELOS_COMP}.
		For instance, we can define the homomorphism $\varphi: R[t] \rightarrow R[z,z^{-1}]$ given by $\varphi(x_i)=x_iz$, $\varphi(y_i)=y_iz$ and $\varphi(t)=1/z^2$, then the restriction ${\varphi\mid}_{\Rees(I)}$ of $\varphi$ to $\Rees(I)$ will give us the required isomorphism because both algebras are integral domains of the same dimension (see \autoref{properties_canonical_module} $(i)$).

		Hence we will identify the ideal $\EEQ$ with the kernel of $\pi$.
		Let $N$ be the incidence matrix of the cone graph $C(G)$. 
		From \cite[Proposition 4.2]{VILLARREAL_EDGE}, we have that a vector $\alpha \in \Ker(N) \cap \mathbb{Z}^{m+n+q}$ is a circuit of $N$ if and only if the monomial walk defined by $\alpha$ corresponds to an even cycle or to two edge disjoint odd cycles joined by a path.
		
		Since the graph $G$ is bipartite, then an odd cycle in $C(G)$ will necessarily contain the vertex $z$.
		Therefore the monomial walks defined by the circuits of $N$ are of the following types:
		\begin{enumerate}[(i)]
			\item An even cycle in $C(G)$ that does not contain the vertex $z$.
			\item An even cycle in $C(G)$ that contains the vertex $z$.
			\item Two odd cycles in $C(G)$ whose intersection is exactly the vertex $z$.
		\end{enumerate}		
	The figure below shows how the cases $(ii)$ and $(iii)$ may look.
	\vspace*{-.1cm}	
	\begin{center}
	\definecolor{xdxdff}{rgb}{0.49019607843137253,0.49019607843137253,1.}
\definecolor{ffttcc}{rgb}{1.,0.2,0.8}
\begin{tikzpicture}[line cap=round,line join=round,>=triangle 45,x=1.0cm,y=1.0cm]
\clip(4.,0.2) rectangle (15.,3.3);
\draw [line width=1.2pt] (5.66,1.99) circle (0.9808159868191381cm);
\draw [line width=1.2pt] (11.64,1.99) circle (0.9808159868191381cm);
\draw [line width=1.2pt] (8.2,1.99) circle (0.9808159868191373cm);
\draw [line width=1.2pt] (13.6,1.95) circle (0.9808159868191381cm);
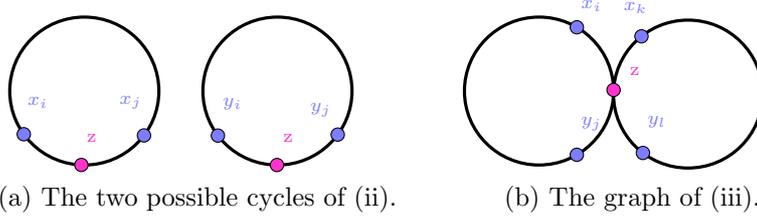
\draw (4.4,0.86) node[anchor=north west] {(a) The two possible cycles of (ii).};
\draw (11.06,0.86) node[anchor=north west] {(b) The graph of (iii).};
\begin{scriptsize}
\draw [fill=ffttcc] (5.621186229764297,1.0099523015484961) circle (2.5pt);
\draw[color=ffttcc] (5.76,1.37) node {z};
\draw [fill=xdxdff] (6.4446527894553105,1.4015104079085172) circle (2.5pt);
\draw[color=xdxdff] (6.27,1.84) node {$x_j$};
\draw [fill=xdxdff] (4.863850028531016,1.4171603863820732) circle (2.5pt);
\draw[color=xdxdff] (5.05,1.84) node {$x_i$};
\draw [fill=ffttcc] (8.19355758292411,1.009205171678491) circle (2.5pt);
\draw[color=ffttcc] (8.34,1.37) node {z};
\draw [fill=xdxdff] (8.992749828595295,1.4124641056504563) circle (2.5pt);
\draw[color=xdxdff] (8.77,1.74) node {$y_j$};
\draw [fill=xdxdff] (7.416978275315728,1.399341910516357) circle (2.5pt);
\draw[color=xdxdff] (7.61,1.82) node {$y_i$};
\draw [fill=ffttcc] (12.620706812337366,2.004633804531277) circle (2.5pt);
\draw[color=ffttcc] (12.9,2.25) node {z};
\draw [fill=xdxdff] (12.137293445772439,2.8353988578131473) circle (2.5pt);
\draw[color=xdxdff] (12.33,3.12) node {$x_i$};
\draw [fill=xdxdff] (12.137293445772439,1.1446011421868525) circle (2.5pt);
\draw[color=xdxdff] (12.33,1.56) node {$y_j$};
\draw [fill=xdxdff] (13.006775374966306,1.1689209103723022) circle (2.5pt);
\draw[color=xdxdff] (13.19,1.58) node {$y_l$};
\draw [fill=xdxdff] (12.98728911047188,2.7158886119101515) circle (2.5pt);
\draw[color=xdxdff] (12.91,3.1) node {$x_k$};
\end{scriptsize}
\end{tikzpicture}
	\end{center}
	
	Since the circuits of the matrices $M$ and $N$ coincide, now we translate these monomial walks in $C(G)$ into binomials of $\EEQ$. 
	An even cycle in $C(G)$ not containing $z$, is also  an even cycle in $G$, and it determines a binomial in $\EEQ$ using \autoref{notation_trans}.
	In the cases $(ii)$ and $(iii)$, we delete vertex $z$ in order to get a subgraph $H$ of $G$. 
	Thus we have that $H$ is either an even path or two disjoint odd paths, and we translate these into binomials in $\EEQ$ using \autoref{notation_trans}.
	\end{proof}
\end{theorem}

\begin{remark}
	Alternatively in \autoref{thm_univ_grobner_basis}, we can see that the matrices $M$ and $N$ have the same kernel because they are equivalent.
	We multiply the last row of $M$ by $-2$ and then we successively add the rows $1,\dots,n+m$ to the last row; with these elementary row operations we transform $M$ into $N$.
\end{remark}

\begin{example}
	Using \autoref{thm_univ_grobner_basis}, the universal Gr\"obner basis of the defining equations of the Rees algebra of the graph in \autoref{examp_graph_two_odd_cycles} is given by
	\begin{align*}
	\big\{x_2y_2T_1-x_1y_1T_2,\; x_2y_3T_1T_3-x_1y_1T_2T_4,\;
	x_3T_2-x_2T_3,\; x_3y_2T_1-x_1y_1T_3,\\ x_3y_3T_1-x_1y_1T_4,\;
	y_3T_3-y_2T_4,\; x_3y_3T_2-x_2y_2T_4\big\}.
		\end{align*}
	It can also be checked in \cite{M2} using the command {\normalfont\fontfamily{cmtt}\selectfont universalGroebnerBasis}.
\end{example}

\begin{corollary}
	Let $G$ be a bipartite graph and $I=I(G)$ be its edge ideal.
	The universal Gr\"obner basis $\mathcal{U}$ of $\EEQ$ consists of square free binomials with degree at most linear in the variables $x_i$'s and at most linear in the variables $y_i$'s.
\end{corollary}

\section{Upper bound for the $xy$-regularity}\label{section3}

In this section we get an upper bound for the $xy$-regularity of $\Rees(I)$, and the important point is that we will choose a special monomial order.
Using the $xy$-regularity we can find an upper bound for the regularity of all the powers of the edge ideal $I$. 

Since most of the upper bounds for the regularity of the powers of edge ideals are based on the technique of even-connection \cite{REG_POWERS_EDGE_BANERJEE}, then a strong motivation for this section is trying to give new tools for the challenging conjecture:

\begin{conjecture}[Alilooee, Banerjee, Beyarslan and H\`a]
	\label{conjecture_powers_bound}
	Let $G$ be an arbitrary graph then 
	$$
	\reg(I(G)^s) \le 2s + \reg(I(G))-2
	$$
	for all $s\ge 1$.
\end{conjecture}

The following theorem will be crucial in our treatment.

\begin{theorem}
	\label{thm_bound_reg_all_powers}
	{\normalfont(\cite[Theorem 5.3]{ROMER_BIGRAD}, \cite[Proposition 10.1.6]{HERZOG_HIBI_BOOK})} The regularity of each power $I^s$ is bounded by 
	$$
	\reg(I^s) \le 2s + \reg_{xy}(\Rees(I)).
	$$
\end{theorem} 

By fixing a particular monomial order $<$ in $S$, then we can see the initial ideal $\iniTerm_{<}(\EEQ)$ as the special fibre of a flat family whose general fibre is $\EEQ$ (see e.g. \cite[Chapter 3]{HERZOG_HIBI_BOOK} or \cite[Chapter 15]{EISENBUD_COMM}), and we can get a bigraded version of \cite[Theorem 3.3.4, (c)]{HERZOG_HIBI_BOOK}.

\begin{theorem}
	\label{thm_bound_init_ideal}
	Let $<$ be a monomial order in $S$, then we have 
	$$
	\reg_{xy}(\Rees(I)) \le \reg_{xy}(S/\iniTerm_{<}(\EEQ)).
	$$
\end{theorem}

Let $\mathcal{M}$ be an arbitrary maximal matching in $G$ with $\lvert \match \rvert = r$. 
We assume that the vertices of $G$ are numbered in such a way that $\match$ consists of the edges 
$$
\match =\big\{\{x_1,y_1\},  \{x_2,y_2\}, \ldots, \{x_r,y_r\}\big\},
$$
and also we assume that $n=\lvert X \rvert \le \lvert Y \rvert = m$.

In $R=\kk[x_1\ldots,x_n,y_1,\ldots,y_m]$ we consider the lexicographic monomial order induced by
$$
x_n>\ldots>x_{2}> x_1 >  y_m>\ldots>y_{2}> y_1.
$$
We choose an arbitrary monomial order $<^{\#}$ on $\kk[T_1,\ldots,T_q]$,
then we define the following monomial order $<^{\match}$ on $S=\kk[x_1,\ldots,x_n,y_1,\ldots,y_m,T_1,\ldots,T_q]$: for two monomials $\mathbf{x}^{\alpha_1}\mathbf{y}^{\beta_1}\mathbf{T}^{\gamma_1}$ and  $\mathbf{x}^{\alpha_2}\mathbf{y}^{\beta_2}\mathbf{T}^{\gamma_2}$ we have 
$$
\mathbf{x}^{\alpha_1}\mathbf{y}^{\beta_1}\mathbf{T}^{\gamma_1} <^{\match} \mathbf{x}^{\alpha_2}\mathbf{y}^{\beta_2}\mathbf{T}^{\gamma_2}
$$
if either 
\begin{enumerate}[(i)]
	\item $\mathbf{x}^{\alpha_1}\mathbf{y}^{\beta_1} < \mathbf{x}^{\alpha_2}\mathbf{y}^{\beta_2}$ or
	\item $\mathbf{x}^{\alpha_1}\mathbf{y}^{\beta_1} = \mathbf{x}^{\alpha_2}\mathbf{y}^{\beta_2}$ and $\mathbf{T}^{\gamma_1} <^{\#}\mathbf{T}^{\gamma_2}$.
\end{enumerate}

Let $\GB$ be the reduced Gr\"obner basis of $\EEQ$ with respect to $<^\match$. 
The possible type of binomials inside $\GB$ were described in \autoref{thm_univ_grobner_basis},  now we focus on obtaining a more refined information about the type $(iii)$ in \autoref{notation_trans}.

\begin{notation}
	\label{nota_specif_binomial}
	In this section, for notational purposes (and without loss of generality) we shall assume that $w_1$ and $w_2$ are disjoint odd paths of the form
	\begin{align*}
	w_1&=(x_e,u_1,\ldots,u_{2a},y_f),\\
	w_2&=(x_g,v_1,\ldots,v_{2b},y_h).
	\end{align*}
	Then we analyze the binomial $x_ey_fT_{(w_1,w_2)^+}- x_gy_hT_{(w_1,w_2)^-}$.	
\end{notation}

\begin{lemma}
	\label{lem_one_end_point}
	Let $x_ey_fT_{(w_1,w_2)^+}- x_gy_hT_{(w_1,w_2)^-} \in \GB$, then we have
	\begin{enumerate}[(i)]
		\item at least one of the vertices $x_e, y_f$ is in the matching $\match$, i.e. $e \le r$ or $f \le r$; 
		\item at least one of the vertices $x_g, y_h$ is in the matching $\match$, i.e. $g \le r$ or $h \le r$.
	\end{enumerate}
	\begin{proof}
		$(i)$ First, assume that $a=0$, i.e. $w_1$ has length one. 
		Since $\match$ is a maximal matching then we necessarily get that $e \le r$ or $f \le r$.
		
		Now let $a>0$, and by contradiction assume that $e > r$ and $f > r$.
		From the maximality of $\match$, we get that $u_1=y_j$ where $j \le r$.
		We consider the even path 
		$$
		w_3=(y_j, \ldots, u_{2a}, y_f),
		$$
		then using \autoref{notation_trans} we get the binomial 
		$$
			F=y_jT_{w_3^+} - y_fT_{w_3^-} \in \EEQ.
		$$
		We have $\iniTerm_{<^\match}(F)=y_fT_{w_3^-}$ because $f> j$.
		So we obtain that $\iniTerm_{<^\match}(F)$ divides $x_ey_fT_{(w_1,w_2)^+}$, and this contradicts that $\GB$ is reduced.
		
		$(ii)$ Follows identically.
	\end{proof}
\end{lemma}

In the rest of this note we assume the following.

\begin{notation}
	$b(G)$ represents the minimum cardinality of the maximal matchings of $G$ and $\normalfont\text{match}(G)$ denotes the maximum cardinality of the matchings of $G$
\end{notation}

\begin{theorem}
	\label{thm_bound_xy_reg}
	Let $G$ be a bipartite graph and $I=I(G)$ be its edge ideal.
	The $xy$-regularity of $\Rees(I)$ is bounded by 
	$$
	\reg_{xy}(\Rees(I)) \le \min\big\{\lvert X \rvert - 1, \; \lvert Y \rvert - 1, \; 2b(G)-1 \big \}.
	$$
	\begin{proof}
		From \autoref{thm_bound_init_ideal}, it is enough to prove that
		 $$
		 \reg_{xy}(S/\iniTerm_{<^\match}(\EEQ))  \le \min\big\{\lvert X \rvert - 1, \; \lvert Y \rvert - 1, \; 2r-1 \big \}.
		 $$
		Let $\{m_1,\ldots,m_c \}$ be the monomials obtained as the initial terms of the elements of $\GB$.
		We consider the Taylor resolution (see e.g. \cite[Section 7.1]{HERZOG_HIBI_BOOK})
		$$
		0 \longrightarrow T_c \longrightarrow \cdots \longrightarrow T_1 \longrightarrow T_0 \longrightarrow S/\iniTerm_{<^\match}(\EEQ) \longrightarrow 0,
		$$
		where each $T_i$ as a bigraded $S$-module has the structure
		$$
		T_i=\bigoplus_{1 \le j_1 < \ldots < j_i \le c} S\big(-\deg_{xy}(\text{lcm}(m_{j_1},\ldots,m_{j_i})), -\deg_{T}(\text{lcm}(m_{j_1},\ldots,m_{j_i})) \big).
		$$
		
		From it, we get the upper bound
		$$
		\reg_{xy}(S/\iniTerm_{<^\match}(\EEQ)) \le \max\big\{\deg_{xy}(\text{lcm}(m_{j_1},\ldots,m_{j_i})) - i \mid \{j_1,\ldots,j_i\} \subset \{1,\ldots,c \} \big\}.
		$$
		When $\deg_{xy}(m_{j_i}) \le 1$, then we have 
		\begin{equation}
			\label{reduction_degree_xy}
			\deg_{xy}(\text{lcm}(m_{j_1},\ldots,m_{j_i})) - i \le 	\deg_{xy}(\text{lcm}(m_{j_1},\ldots,m_{j_{i-1}})) - (i-1).
		\end{equation}
		
		So, according with \autoref{thm_univ_grobner_basis}, we only need to consider subsets $\{j_1,\ldots,j_i\}$ such that for each $1 \le k \le i$ we have  $m_{j_k}=\iniTerm_{<^\match}(F_k)$ and $F_k$ is a binomial
		as in \autoref{nota_specif_binomial}.
		We use the notation $\iniTerm_{<^\match}(F_k)=x_{e_k}y_{f_k}B_k$,  where $B_k$ is a monomial in the $T_i$'s.
		Also, we can assume that $x_{e_1}y_{f_1}, x_{e_2}y_{f_2}, \ldots, x_{e_k}y_{f_k}$ are pairwise relatively prime, 
		because we can make a reduction like in \autoref{reduction_degree_xy} if this condition is not satisfied.
		
		Thus, in order to finish the proof, we only need to show that we necessarily have $i \le \min\{\lvert X \rvert - 1, \; \lvert Y \rvert - 1, \; 2r-1 \}$ under the two  previous conditions.
		Since the two paths that define each $F_k$ are disjoint, then by the monomial order chosen we have that  $e_k > 1$ for each $k$, and by a ``pigeonhole'' argument follows that $i \le \lvert X \rvert - 1 \le \lvert Y \rvert- 1$.
		Also, from \autoref{lem_one_end_point} there are at most $2r-1$ available positions to satisfy the condition of being co-primes.
		Thus we have $i \le 2r-1$, and the result of the theorem follows because $\match$ is an arbitrary maximal matching.
		\end{proof}
\end{theorem}

\begin{corollary}
	\label{cor_upper_bound_bipart_order}
	Let $G$ be a bipartite graph and $I=I(G)$ be its edge ideal.
	For all $s \ge 1$ we have 
	$$
	\reg(I^s) \le 2s + \min\big\{\lvert X \rvert - 1, \; \lvert Y \rvert - 1, \; 2b(G)-1 \big \}.
	$$
	\begin{proof}
		It follows from \autoref{thm_bound_xy_reg} and \autoref{thm_bound_reg_all_powers}.		
	\end{proof}
\end{corollary}

\begin{remark}
	\label{comparison}
	From the fact that $\normalfont\text{co-chord}(G) \le \text{match}(G) \le \min\{\lvert X\rvert, \lvert Y \rvert\}$ (see \cite{REG_POWERS_BIPART}) and $\normalfont\text{match}(G) \le 2b(G)$ (see \cite[Proposition 2.1]{HIBI_MATCHING}), then we have the following relations
	\begin{equation*}
	\normalfont\text{co-chord}(G) - 1\le \text{match}(G) -1 \le \min\big\{\lvert X \rvert - 1, \; \lvert Y \rvert - 1, \; 2b(G)-1 \big \}.
	\end{equation*}
\end{remark}

Although the last upper bound is weaker, it is interesting that an approach based on Gr\"obner bases can give a sharp answer in several cases.


In the last part of this section we deal with the case of a complete bipartite graph.
The Rees algebra of these graphs was studied in \cite{VILLARREAL_COMPLETE_BIPART}.

\begin{notation}
	\label{nota_complete_bipartite}
	By $G$ we will denote a complete bipartite graph with bipartition $X=\{x_1,\ldots,x_n\}$ and $Y=\{y_1,\ldots,y_m\}$.
	Let $I=\{x_iy_j \mid 1\le i \le n,1 \le j \le m\}$ be the edge ideal of $G$ and let $T_{ij}$ be the variable that corresponds to the edge $x_iy_j$. Thus we have a canonical map
	\begin{align}
	\label{def_complete_Rees}
	\begin{split}
	S&=\kk[x_i\text{'s},\;y_j\text{'s},\; T_{ij}\text{'s}] \xrightarrow{\psi} \Rees(I) \subset R[t],\\ \quad
	&\psi(x_i) = x_i,\quad \psi(y_i) = y_i, \quad \psi(T_{ij}) = x_iy_jt.
	\end{split}
	\end{align}
	Let $\EEQ$ be the kernel of this map.
	For simplicity of notation we keep the same monomial order $<^{\match}$.
\end{notation}

Exploiting our characterization of the universal Gr\"obner basis of $\EEQ$, 
we shall prove that all the powers of the edge ideal of $G$ have a linear free resolution.

\begin{lemma}
	\label{grob_basis_complete}
	Let $G$ be a complete bipartite graph.
	 The reduced Gr\"obner basis $\GB$ consists of binomials with linear $xy$-degree.
	\begin{proof}
		From \autoref{thm_univ_grobner_basis} we only need to show that any binomial determined by two disjoint odd paths is not contained in $\GB$.
		Let $x_ey_fT_{(w_1,w_2)^+}- x_gy_hT_{(w_1,w_2)^-}$ be a binomial like in \autoref{nota_specif_binomial}. 
		By contradiction assume that $x_ey_fT_{(w_1,w_2)^+}- x_gy_hT_{(w_1,w_2)^-} \in \GB$.
		
		Without loss of generality we assume that $e > g$.
		Since $G$ is complete bipartite, we choose the edge $x_ey_h$ and we append it to $w_2$, that is
		$$
		w_3 = (x_g,v_1,\ldots,v_{2b},y_h, x_e).
		$$
		Using \autoref{notation_trans} we get the binomial 
		$$
			F=x_gT_{w_3^+} - x_eT_{w_3^-} \in \EEQ,
		$$
		with initial term $\iniTerm_{<^\match}(F)=x_eT_{w_3^-}$ because $e > g$.
		Thus we get that $\iniTerm_{<^\match}(F)$ divides $x_ey_fT_{(w_1,w_2)^+}$, a contradiction.						
	\end{proof}
\end{lemma}

\begin{corollary}
	\label{reg_complete_bipart}	
	Let $G$ be a complete bipartite graph and $I=I(G)$ be its edge ideal.  For all $s\ge 1$ we have $\reg(I^s)=2s$.
	\begin{proof}
	Using \autoref{grob_basis_complete} and repeating the same argument of \autoref{thm_bound_xy_reg} we can get $\reg_{xy}(\Rees(I))=0$.
	Again, the result follows by \autoref{thm_bound_reg_all_powers}.
	\end{proof}
\end{corollary}

We remark that this previous result also follows from \cite{REG_POWERS_BIPART} since it is easy to check that $\text{co-chord}(G)=1$ (i.e. it is a co-chordal graph) in the case of complete bipartite graphs.

\section{The total regularity of $\Rees(I)$}\label{section4}


In the previous sections we heavily exploited the fact that the matrix $M$ (corresponding to $\Rees(I)$) is totally unimodular in the case of a bipartite graph $G$.
From \cite[Theorem 2.1]{VILLARREAL_BIPARTITE} we have that $\Rees(I)$ is a normal domain, then a famous theorem by Hochster \cite{HOCHSTER} (see e.g. \cite[Theorem 6.10]{BRUNS_GUB_POLYTOPES} or \cite[Theorem 6.3.5]{BRUNS_HERZ}) implies that $\Rees(I)$ is Cohen-Macaulay.
So, the Rees algebra $\Rees(I)$ of a bipartite graph $G$ is also special from a more algebraic point of view (see \cite{SIMIS_VASC_VILLARREAL_IDEAL_GRAPH}).


For notational purposes we let $N$ be $N=n+m$.
It is well known that the canonical module of $S$ (with respect to our bigrading) is given by $S(-N, -q)$ (see e.g. \cite[Proposition 6.26]{BRUNS_GUB_POLYTOPES}, or \cite[Example 3.6.10]{BRUNS_HERZ} in the $\mathbb{Z}$-graded case). 
The Rees cone  is the polyhedral cone of $\mathbb{R}^{N+1}$ generated by the set of vectors
$$
\mathcal{A} = \big\{ v \mid v \text{ is a column of } M \text{ in \autoref{matrix_def_Rees}} \big\},
$$
and we will denote it by $\mathbb{R}_+\mathcal{A}$.
The irreducible representation of the Rees cone for a bipartite graph was given in \cite[Section 4]{VILLARREAL_BIPARTITE}. 

\begin{proposition}
	\label{properties_canonical_module}
		Adopt \autoref{general_notation}.
	The following statements hold:
	\begin{enumerate}[(i)]
		\item The Krull dimension of $\Rees(I)$ is $\dim(\Rees(I))=N+1$.
		\item The projective dimension of $\Rees(I)$ as an $S$-module is equal to the number of edges minus one, that is, $p={\normalfont\text{pd}}_S(\Rees(I))=q-1$.
		\item The canonical module of $\Rees(I)$ is given by 
		$$
		\omega_{\Rees(I)} = {}^*\Ext_S^p\big(\Rees(I),\; S(-N,-q)\big).
		$$
		\item The bigraded Betti numbers of $\Rees(I)$ and $\omega_{\Rees(I)}$ are related by
		$$
		\beta_{i, (a, b)}^S(\Rees(I)) = \beta_{p-i, (N-a, q-b)}^S(\omega_{\Rees(I)}).
		$$
	\end{enumerate}
	\begin{proof}
		$(i)$ The Rees cone $\mathbb{R}_+\mathcal{A}$ has dimension $N+1$ and the Krull dimension of $\Rees(I)$ is equal to it (see e.g. \cite[Lemma 4.2]{STURMFELS_MONOG}). 
		More generally, it also follows from \cite[Proposition 2.2]{SIMIS_ULRICH_VASC_REES_MOD}.

		Since clearly $\Rees(I)$ is a finitely generated $S$-module, then the statements $(ii)$ and $(iii)$ follow from \cite[Theorem 6.28]{BRUNS_GUB_POLYTOPES} (see \cite[Proposition 3.6.12]{BRUNS_HERZ} for the $\mathbb{Z}$-graded case).
		
		The statement $(iv)$ follows from \cite[Theorem 6.18]{BRUNS_GUB_POLYTOPES}; also, see \cite[page 224, equation 6.6]{BRUNS_GUB_POLYTOPES}.
	\end{proof}
\end{proposition}

	Due to a formula of Danilov and Stanley (see e.g. \cite[Theorem 6.31]{BRUNS_GUB_POLYTOPES} or \cite[Theorem 6.3.5]{BRUNS_HERZ}), the canonical module of $\Rees(I)$ is the ideal given by 
$$
\omega_{\Rees(I)}=\big(\{ x_1^{a_1}\cdots x_n^{a_n}y_1^{a_{n+1}}\cdots y_N^{a_N} t^{a_{N+1}} \mid a = (a_i) \in (\mathbb{R}_+\mathcal{A})^{\circ} \cap \mathbb{Z}^{N+1} \}\big),
$$		
where $(\mathbb{R}_+\mathcal{A})^{\circ}$ denotes the topological interior of $\mathbb{R}_+\mathcal{A}$.

Now we can compute the total regularity of $\Rees(I)$.

\begin{theorem}
	\label{reg_Rees}
		Let $G$ be a bipartite graph and $I=I(G)$ be its edge ideal.
	The total regularity of $\Rees(I)$ is given by 
	$$
	\reg(\Rees(I)) = \text{match}(G).
	$$
	\begin{proof}		
		In the case of the total regularity, we can see $\Rees(I)$ as a standard graded $S$-module (i.e. $\deg(x_i)=\deg(y_i)=\deg(T_i)=1$), and since $\Rees(I)$ is a Cohen-Macaulay $S$-module then the regularity can be computed with the last Betti numbers (see e.g. \cite[page 283]{SCHENZEL_NOTES} or \cite[Exercise 20.19]{EISENBUD_COMM}). 
		Thus, from 	\autoref{properties_canonical_module} we get 
		\begin{align*}
			\reg(\Rees(I)) &= \max\big\{ a+b-p \mid \beta_{p, (a, b)}^S(\Rees(I))\neq 0
			\big\}\\
			&= \max\big\{ a+b-p\mid \beta_{0, (N-a, q-b)}^S(\omega_{\Rees(I)}) \neq 0 \big\}\\
			&= N +1 - \min\big\{ a+b \mid \beta_{0, (a, b)}^S(\omega_{\Rees(I)})\neq 0  \big\},
		\end{align*}
		and by the bigrading that we are using ($\bideg(x_i)=\bideg(y_i)=(1,0)$ and $\bideg(t)=(-2,1)$) then we obtain 
		$$
		\reg(\Rees(I)) = N+1-\min\big\{a_1+\cdots+a_N-a_{N+1}  \mid a = (a_i) \in (\mathbb{R}_+\mathcal{A})^{\circ} \cap \mathbb{Z}^{N+1} \big\}.
		$$
		One can check that the number 
		$$
		-\min\big\{a_1+\cdots+a_N-a_{N+1}  \mid a = (a_i) \in (\mathbb{R}_+\mathcal{A})^{\circ} \cap \mathbb{Z}^{N+1} \big\}
		$$
		coincides with the $a$-invariant of $\Rees(I)$ with respect to the $\mathbb{Z}$-grading induced by $\deg(x_i)=\deg(y_i)=1$ and $\deg(t)=-1$. 
		This last formula can be evaluated with the irreducible representation of the Rees cone \cite[Corollary 4.3]{VILLARREAL_BIPARTITE}, it was done in \cite[Proposition 4.5]{VILLARREAL_BIPARTITE}, and from it we get 
		$$
		\reg(\Rees(I)) = N - \beta_0, 
		$$
		where $\beta_0$ denotes the maximal size of an independent set of $G$.
		The minimal size of a vertex cover is equal to $N-\beta_0$, and we finally get
		$$
		\reg(\Rees(I)) = \text{match}(G)
		$$
		from K\"onig's theorem.
	\end{proof}
\end{theorem}

The following bound was obtained for the first power of the edge ideal in \cite[Theorem 6.7]{TAI_ADAM_HYPERGRAPHS}.
\begin{corollary}
	\label{cor_upper_bound_all_powers}
	Let $G$ be a bipartite graph and $I=I(G)$ be its edge ideal.
	For all $s \ge 1$ we have 
	$$
	\reg(I^s) \le 2s + \text{match}(G) -1.
	$$
	\begin{proof}
		It is enough to prove that $\reg_{xy}(\Rees(I)) \le \reg(\Rees(I)) - 1$.
		In the minimal bigraded free resolution \autoref{minimal_free_res_Rees} of $\Rees(I)$, suppose that $\reg_{xy}(\Rees)=a_{ij}-i$ for some $i,j \in \mathbb{N}$.
		Since necessarily $b_{ij} \ge 1$ and 
		$$
		a_{ij}+b_{ij}-i \le \reg(\Rees(I)),
		$$
		then we get the expected inequality.
	\end{proof}
\end{corollary}

This previous upper bound is sharp in some cases (see \cite[Lemma 4.4]{FORESTS_CYCLES}).
In the following corollary we get information about the eventual linearity.

\begin{corollary}
	\label{cor_eventual_linearity}
		Let $G$ be a bipartite graph and $I=I(G)$ be its edge ideal.
	For all $s \ge {\normalfont\text{match}}(G)+q+1$ we have
	$$
	\reg(I^{s+1})=\reg(I^s)+2.
	$$
	\begin{proof}
		With the same argument of \autoref{cor_upper_bound_all_powers} we can prove that $\reg_{T}(\Rees(I)) \le \reg(\Rees(I))$, here the difference is that in the minimal bigraded free resolution \autoref{minimal_free_res_Rees} we can have free modules of the type $S(0,-b_{ij})$ (for instance, in the syzygies of $\Rees(I)$ the ones that come from even cycles).
		Then the statement of the corollary follows from \cite[Proposition 3.7]{CUTKOSKY_HERZOG_TRUN_LINEARITY_REG}.
	\end{proof}
\end{corollary}

\section{Some final thoughts}
\label{section5}

In the last part of this note we give some ideas and digressions about \autoref{conjecture_powers_bound}.
Using a ``refined Rees approach'' with respect to the one of this note, one might get an answer to this conjecture for general graphs or perhaps for special families of graphs:
\begin{itemize}
	\item Restricting the minimal bigraded free resolution \autoref{minimal_free_res_Rees} of $\Rees(I)$ to a graded $T$-part gives an exact sequence
	$$
	0 \longrightarrow (F_p)_{(*,k)} \longrightarrow \cdots \longrightarrow (F_1)_{(*,k)} \longrightarrow (F_0)_{(*,k)} \longrightarrow (\Rees(I))_{(*,k)} \longrightarrow 0
	$$
	for all $k$.
	This gives a (possibly non-minimal)  graded free $R$-resolution of 
	$$(\Rees(I))_{(*,k)} \cong I^k(2k).
	$$
	But in the case $k=1$ one can check that 
	$$
	0 \longrightarrow (F_p)_{(*,1)} \longrightarrow \cdots \longrightarrow (F_1)_{(*,1)} \longrightarrow (F_0)_{(*,1)} \longrightarrow I(2) \longrightarrow 0
	$$
	is indeed the minimal free resolution of $I(2)$.
	Thus, one can read the regularity $I$ from \autoref{minimal_free_res_Rees}, and a solution to \autoref{conjecture_powers_bound} can be given by proving that 
	$$
	\max_{i,j}\big\{ a_{ij}-i \big\} = \max_{i,j}\big\{ a_{ij}-i \mid b_{ij}=1\big\}.
	$$	
	\item For bipartite graphs, Gr\"obner bases techniques can give very good results (for instance, in the case of complete bipartite graphs). Perhaps, for special families of bipartite graphs one can give ``good'' monomial orders.  
	
	\item The existence of a canonical module in the case of bipartite graphs could give more information about the minimal bigraded free resolution of $\Rees(I)$. 
	From \cite[Theorem 7.26]{BRUNS_GUB_POLYTOPES} we have that the maximal $xy$-degree and the maximal $T$-degree on each $F_i$ of \autoref{minimal_free_res_Rees} form  weakly increasing sequences of integers, that is 
	$$
	\max_j\{a_{ij}\} \le \max_j\{a_{i+1,j}\} \quad \text{ and } \quad 	\max_j\{b_{ij}\} \le \max_j\{b_{i+1,j}\}
	$$
	(see e.g. \cite[Exercise 20.19]{EISENBUD_COMM} for the $\mathbb{Z}$-graded case).
	Thus a more detailed analysis of the polyhedral geometry of the Rees cone $\mathbb{R}_+\mathcal{A}$ could give better results.
\end{itemize}

\section*{Acknowledgments}
The work presented in this note started thanks to the PRAGMATIC 2017 Research School in Algebraic Geometry and Commutative Algebra ``Powers of ideals and ideals of powers'' held in Catania, Italy, in June 2017.
The author is very grateful to the organizers of the event, and to the professors Brian Harbourne, Adam Van Tuyl, Enrico Carlini and T\`ai Huy H\`a who gave the lectures. 
The author is specially grateful to T\`ai Huy H\`a for his support and for insisting on a ``Rees approach''.
The author is grateful to Carlos D'Andrea and Aron Simis for useful suggestions.
The use of \textit{Macaulay2} \cite{M2} was very important in the preparation of this note.
The author wishes to thank the referee for numerous suggestions to improve the exposition.

\addcontentsline{toc}{section}{Bibliography}
\bibliographystyle{elsarticle-num} 


\end{document}